\documentstyle [12pt,twoside, epsf]{article}

\def\picill#1by#2(#3)
{\vbox to #2
{\hrule width #1 height 0pt depth 0pt
\vfill\epsffile{#3}}}

\begin{document}

\date{}

\title{\Large\bf Virtual Knot Theory}

\author{Louis H. Kauffman \\ Department of Mathematics, Statistics and Computer
Science \\ University of Illinois at Chicago \\ 851 South Morgan Street\\
Chicago, IL, 60607-7045}

\maketitle

\thispagestyle{empty}

\subsection*{\centering Abstract}

{\em This paper is an introduction to the theory of virtual knots. It is
dedicated to the memory of Francois Jaeger. }

\section{Introduction}

This paper is an introduction to the subject of virtual knot theory, a
generalisation of classical knot theory that I discovered in 1996
\cite{Lectures}. This paper gives the basic definitions, some fundamental
properties and  a collection of examples. Subsequent papers will treat specific
topics such as classical and quantum link invariants and Vassiliev invariants for
virtual knots and links in more detail. \vspace{3mm}

Throughout this paper I shall refer to knots and links by the generic term
``knot".  In referring to a trivial fundamental group of a knot, I mean that the
fundamental group is isomorphic to the integers.

\vspace{3mm}

The paper is organised as follows.  Section 2 gives the definition of a virtual
knot in terms of diagrams and moves on diagrams. Section 3 discusses both the
motivation from knots in thickened surfaces and the abstract properties of Gauss
codes. Section 3 proves basic results about virtual knots by using reconstruction
properties of Gauss codes. In particular, we show how virtual knots can be
identified as virtual by examining their codes.  Section 4 discusses the
fundamental group and the quandle extended for virtual knots. Examples are given
of non-trivial virtual knots with trivial (isomorphic to the integers)
fundamental group. An example shows that some virtual knots are distinguished
from their mirror images by the fundamental group, a very non-classical effect.
Section 5 shows how the bracket polynomial (hence the Jones polynomial) extends
naturally to virtuals and gives examples of non-trivial virtual knots with
trivial Jones polynomial. Examples of infinitely many distinct virtuals with the
same fundamental group are verified by using the bracket polynomial. An example
is given of a knotted virtual with trivial fundamental group and unit Jones
polynomial. It is conjectured that this phenomenon cannot happen with virtuals
whose shadow code is classical.  In section 6 we show how to extend quantum link
invariants and introduce the concept of virtual framing.  This yields a virtually
framed bracket polynomial distinct from the model in the previous section and to
generalisation of this model to an invariant, $\overline{Z}(K),$ of virtual
regular isotopy depending on infinitely many variables. Section 7 discusses
Vassiliev invariants, defines graphical finite type and proves that the weight
systems are finite for the virtual Vassiliev invariants arising from the Jones
polynomial. Section 8 is a discussion of open problems. \vspace{3mm}

\noindent {\bf Acknowledgement.}    It gives the author pleasure to thank the
National Science Foundation for support of this research under NSF Grant
DMS-9205277,the NSA for  partial support under grant number MSPF-96G-179 and the
Mittag-Leffler Institute for hospitality during the writing of this paper.
\vspace{3mm}

\section{Defining Virtual Knots and Links}

A classical knot \cite{Ashley} can be represented by a diagram. The diagram is a
4-regular plane graph with extra structure at its nodes. The extra structure is
classically intended to indicate a way to embed a circle in three dimensional
space. The shadow of a projection of this embedding is the given plane graph.
Thus we are all familiar with the usual convention for illustrating a crossing by
omitting a bit of arc at the node of the plane graph. The bit omitted is
understood to pass underneath the uninterrupted arc. See Figure 1 .\vspace{3mm}

From the point of view of a topologist, a knot diagram represents an ``actual"
knotted (possibly unknotted) loop embedded in three space. The crossing structure
is an artifact of the projection to the plane. \vspace{3mm}

I shall define a virtual knot (or link) diagram.  The definition of a virtual
diagram is just this: We allow a new sort of crossing, denoted as shown in Figure
1  as a 4-valent vertex with a small circle around it. \vspace{3mm}

$$ \picill3inby2in(F1.EPSF) $$
\begin{center}
{\bf Figure 1 - Crossings and Virtual Crossings}
\end{center}

%\begin{figure}[htbp] \vspace*{40mm} \special{pntg=F1.pntg} \vspace*{13pt}
%\begin{center} { Figure 1 --- Crossings and Virtual Crossings} \end{center}
%\end{figure} \vspace{3mm}

%\begin{figure}[htbp] \vspace*{140mm} \special{pntg=F2.pntg} \vspace*{13pt}
%\begin{center} { Figure 2 --- Generalised Reidemeister Moves for Virtual Knot
%Theory} \end{center} \end{figure} \vspace{3mm}

$$ \picill3inby6.5in(F2.EPSF) $$
\begin{center}
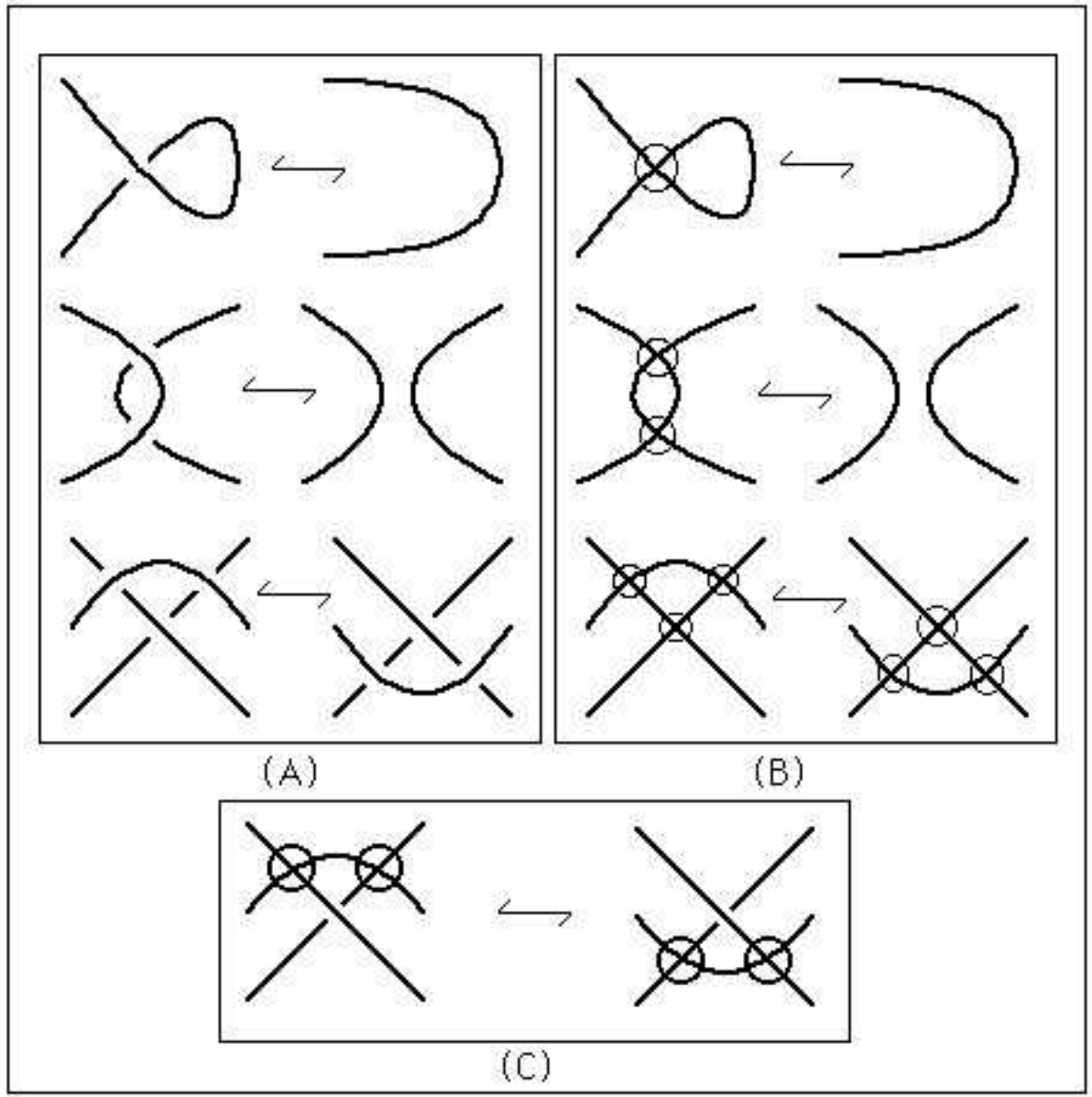
{\bf  Figure 2 --- Generalised Reidemeister Moves for Virtual Knot
Theory}
\end{center}

\noindent This sort of crossing is called virtual. It comes in only one flavor.
You cannot switch over and under in a virtual crossing. However the idea is not
that a virtual crossing is just an ordinary graphical vertex.  Rather, the idea
is that the virtual crossing is not really there. \vspace{3mm}

If I draw a non-planar graph in the plane it necessarily acquires virtual
crossings. These crossings are not part of the structure of the graph itself.
They are artifacts of the drawing of the graph in the plane. The graph theorist
often gets rid of a crossing in the plane by making it into a knot theorist's
crossing, thereby indicating a particular embedding of the graph in three
dimensional space. This is just what we do not do with our  virtual knot
crossings, for then they would be indistinct from  classical crossings.  The
virtual crossings are not there.  We shall make sense of that property by the
following axioms generalising classical Reidemeister moves. See Figure 2.
\vspace{3mm}

The moves fall into three types: (A) Classical Reidemeister moves relating
classical crossings; (B) Shadowed versions of Reidemeister moves relating only
virtual crossings; (C) A triangle move that relates two virtual crossings and one
classical crossing. \vspace{3mm}

\noindent The last move (type C) is the embodiment of our principle that the
virtual crossings are not really there. Suppose that an arc is free of classical
crossings. Then that arc can be arbitrarily moved (holding its endpoints fixed)
to any new location. The new location will reveal a new set of virtual crossings
if the arc that is moved is placed transversally to the remaining part of the
diagram.  See Figure 3 for illustrations of this process and for an example of
unknotting of a virtual diagram. \vspace{3mm}

%\begin{figure}[htbp] \vspace*{140mm} \special{pntg=F3.pntg} \vspace*{13pt}
%\begin{center} { Figure 3 --- Virtual Moves } \end{center} \end{figure}
%\vspace{3mm}

$$ \picill3inby6.5in(F3.EPSF) $$
\begin{center}
{\bf Figure 3 --- Virtual Moves }
\end{center}

The theory of virtual knots is constructed on this combinatorial basis - in terms
of the generalised Reidemeister moves. We will make invariants of virtual knots
by finding functions well-defined on virtual diagrams that are unchanged under
the application of the virtual moves.  The remaining sections of this paper study
many instances of such invariants. \vspace{3mm}

\section{Motivations}

While it is clear that one can make a formal generalisation of knot theory in the
manner so far described, it may not be yet clear why one should generalize in
this particular way.  This section explains two sources of motivation.  The first
is the study of knots in thickened surfaces of higher genus (classical knot
theory is actually the theory of knots in a thickened two-sphere).  The second is
the extension of knot theory to the purely combinatorial domain of Gauss codes
and Gauss diagrams. It is in this second domain that the full force of the
virtual theory comes into play. \vspace{3mm}

\subsection{Surfaces}

Consider the two examples of virtual knots in Figure 4. We shall see later in
this paper that these are both non-trivial knots in the virtual category. In
Figure 4 we have also illustrated how these two diagrams can be drawn (as knot
diagrams) on the surface of a torus. The virtual crossings are then seen as
artifacts of the projection of the torus to the plane. \vspace{3mm}

%\begin{figure}[htbp] \vspace*{140mm} \special{pntg=F4.pntg} \vspace*{13pt}
%\begin{center} { Figure 4 --- Two Virtual Knots } \end{center} \end{figure}
%\vspace{3mm}

$$ \picill3inby5.5in(F4.EPSF) $$
\begin{center}
{\bf  Figure 4 --- Two Virtual Knots}
\end{center}

The knots drawn on the toral surface represent knots in the three manifold  $T
\times I$ where $I$ is the unit interval and $T$ is the torus. If $S_{g}$ is a
surface of genus $g$, then the knot theory in $S_{g} \times I$ is represented by
diagrams drawn on $S_{g}$ taken up to the usual Reidemeister moves transferred to
diagrams on this surface. \vspace{3mm}

As we shall see in the next section, abstract invariants of virtual knots can be
interpreted as invariants for knots that are specifically embedded in $S_{g}
\times I$  for some genus $g.$  The virtual knot theory does not demand the use
of a particular surface embedding, but it does apply to such embeddings. This
constitutes one of the motivations. \vspace{3mm}

\subsection{ Gauss Codes}

A second motivation comes from the use of so-called {\em Gauss codes} to
represent knots and links.  The Gauss code is a sequence of labels for the
crossings with each label repeated twice to indicate a walk along the diagram
from a given starting point and returning to that point.  In the case of multiple
link components we mean a sequence labels, each repeated twice and intersticed by
partition symbols ``/" to indicate the component circuits for the code.
\vspace{3mm}

A {\em shadow} is the projection of a knot or link on the plane with transverse
self-crossings and no information about whether the crossings are overcrossings
or undercrossings. In other words, a shadow is a 4-regular plane graph. On such a
graph we can count circuits  that always cross (i.e., they never use two adjacent
edges in succession at a given vertex) at each crossing that they touch. Such
circuits will be called the {\em components} of the shadow since they correspond
to the components of a link that projects to the shadow. \vspace{3mm}

A single component shadow has a Gauss code that consists in a sequence of
crossing labels, each repeated twice. Thus the trefoil shadow has code  123123. 
A multi-component  shadow has as many sequences as there are components. For
example 12/12 is the code for the Hopf link shadow. \vspace{3mm}

Along with the labels for the crossings one can add the symbols O and U to
indicate that the passage through the crossing was an overcrossing (O) or an
undercrossing (U).  Thus

$$123123$$

\noindent is a Gauss code for the shadow of a trefoil knot and

$$O1U2O3U1O2U3$$

\noindent is a Gauss code for the trefoil knot. \vspace{3mm}

\noindent The Hopf link itself has the code $O1U2/U1O2.$ See Figure 5.
\vspace{3mm}

%\begin{figure}[htbp] \vspace*{140mm} \special{pntg=F5.pntg} \vspace*{13pt}
%\begin{center} { Figure 5 --- Planar and Nonplanar Codes } \end{center}
%\end{figure} \vspace{3mm}

$$ \picill3inby5in(F5.EPSF) $$
\begin{center}
{\bf  Figure 5 --- Planar and Nonplanar Codes}
\end{center}

Suppose that $g$ is such a sequence of labels and that $g$ is free of any
partition labels.  Every label in $g$ appears twice.  The first necessary
criterion for the planarity of $g$ is given by the following definition and
Lemma. \vspace{3mm}

\noindent {\bf Definition.}  A single component Gauss code $g$ is said to be {\em
evenly intersticed} if there is an {\em even} number of labels in between the two
appearances of any label. \vspace{3mm}

\noindent {\bf Lemma 1.}  If $g$ is a single component planar Gauss code, then
$g$ is evenly intersticed. \vspace{3mm}

\noindent {\bf Proof.} This follows directly from the Jordan curve theorem in the
plane.// \vspace{3mm}

\noindent {\bf Example.} The necessary condition for planarity in this Lemma is
not sufficient.  The code $g = 1234534125$  is evenly intersticed but not planar
as is evident from Figure 5. \vspace{3mm}

Non-planar Gauss codes give rise to an infinite collection of virtual knots.
\vspace{3mm}

Local orientations at the crossings give rise to another phenomenon: virtual
knots whose Gauss codes have planar realisations with different local
orientations from their classical counterparts. \vspace{3mm}

By orienting the knot, one can give orientation signs to each crossing relative
to the starting point of the code---using the convention shown in Figure 6.  This
convention designates each oriented crossing with a {\em sign}  of $+1$ or $-1$. 
We say that the crossing has positive sign if the overcrossing line can be turned
through the smaller angle (of the two vertical angles at the crossing)  to
coincide with the direction of the undercrossing line.  The signed code for the
standard trefoil is

$$t = O1+U2+O3+U1+O2+U3+,$$

\noindent while  the signed code for a figure eight knot is

$$f = O1+U2+O3-U4-O2+U1+O4-U3-.$$

\noindent Here we have appended the signs to the corresponding labels in the
code. Thus, crossing number 1 is positive in the figure eight knot, while
crossing number 4 is negative. See Figure 6  for an illustration corresponding to
these codes. \vspace{3mm}

%\begin{figure}[htbp] \vspace*{80mm} \special{pntg=F6.pntg} \vspace*{13pt}
%\begin{center} { Figure 6 --- Signed Gauss Codes} \end{center} \end{figure}
%\vspace{3mm}

$$ \picill3inby4in(F6.EPSF) $$
\begin{center}
{\bf Figure 6 --- Signed Gauss Codes}
\end{center}

Now consider the effect of changing these signs. For example let

$$g =  O1+U2+O3-U1+O2+U3-. $$

\noindent Then $g$ is a signed Gauss code and as Figure 6 illustrates, the
corresponding diagram is forced to have virtual crossings in order to acommodate
the change in orientation. The codes $t$ and $g$ have the same underlying 
(unsigned) Gauss code $O1U2O3U1O2U3,$  but $g$ corresponds to a virtual knot
while t represents the classical trefoil. \vspace{3mm}

Carrying this approach further, we  {\em define} a virtual knot as an equivalence
class of oriented Gauss codes under abstractly defined Reidemeister moves for
these codes---with no mention of virtual crossings. The virtual crossings become
artifacts of a planar representation of the virtual knot.  The move sets of type
B and C for virtuals are diagrammatic rules that make sure that this
representation of the oriented Gauss codes is faithful.  Note, in  particular,
that the move of type C does not alter the Gauss code.  With this point of view
we see that the signed codes are knot theoretic analogues of the set of all
graphs, and that the classical knot (diagrams) are the analogues of the planar
graphs. This is the fundamental combinatorial motivation for our definitions of
virtual knots and their equivalences. \vspace{3mm}

Since it is useful to have a few more facts about the reconstruction of planar
Gauss codes, we conclude this section with a quick review of that subject.
\vspace{3mm}

\subsection{Gauss Codes and Reconstruction}

In this section we recall an algorithm for reconstructing a planar diagram from
its Gauss code.  This algorithm also detects non-planar codes.  We shall see that
for a planar oriented Gauss code, the orientation signs in the code sequence are
determined up to a small number of choices. Such sign sequences will be called
{\em standard} (with the more technical definition to follow). \vspace{3mm}

\noindent We shall prove the following Theorem. \vspace{2mm}

\noindent {\bf Theorem 2.} \noindent If $K$ is a virtual knot whose underlying
Gauss code is planar and whose sign sequence is standard, then $K$ is equivalent
to a classical knot. \vspace{4mm}

The fundamental problem in Gauss codes is to give an algorithm for determining
whether a given code can be realized by a planar shadow. \vspace{3mm}

We will explain the detection and reconstruction algorithms for single component
codes.  The first necessary condition for planarity for a single component code
is that it be evenly intersticed, as we have already remarked in Lemma 1.
\vspace{3mm}

If a code $g$ is planar then a corresponding code for such an evenly paired
Jordan curve can be produced as follows: Let the labels in $g$ be $1,2,..., n.$ 
Starting with $i=1$ {\em reverse the order of labels in between the two
appearances of i.} Do this successively using $i=1,2,...,n.$  Let $g^{*}$ be the
resulting code. \vspace{3mm}

In Figure 7 we see that the crossings of  a planar shadow $E$ can be smoothed to
obtain a single Jordan curve in the plane. This Jordan curve can be seen as a
circle with doubly repeated labels around its circumference so that some labels
are paired by arcs inside the circle, and the remaining labels are paired by arcs
outside the circle. The corresponding code is  $g^{*}$ as defined above.  In this
form of pairing {\em no two pairing arcs intersect one another.} \vspace{3mm}

%\begin{figure}[htbp] \vspace*{140mm} \special{pntg=F7.pntg} \vspace*{13pt}
%\begin{center} { Figure 7 --- If $g$ is Planar, then $g^{*}$ is Dually Paired. }
%\end{center} \end{figure} \vspace{3mm}

$$ \picill3inby5in(F7.EPSF) $$
\begin{center}
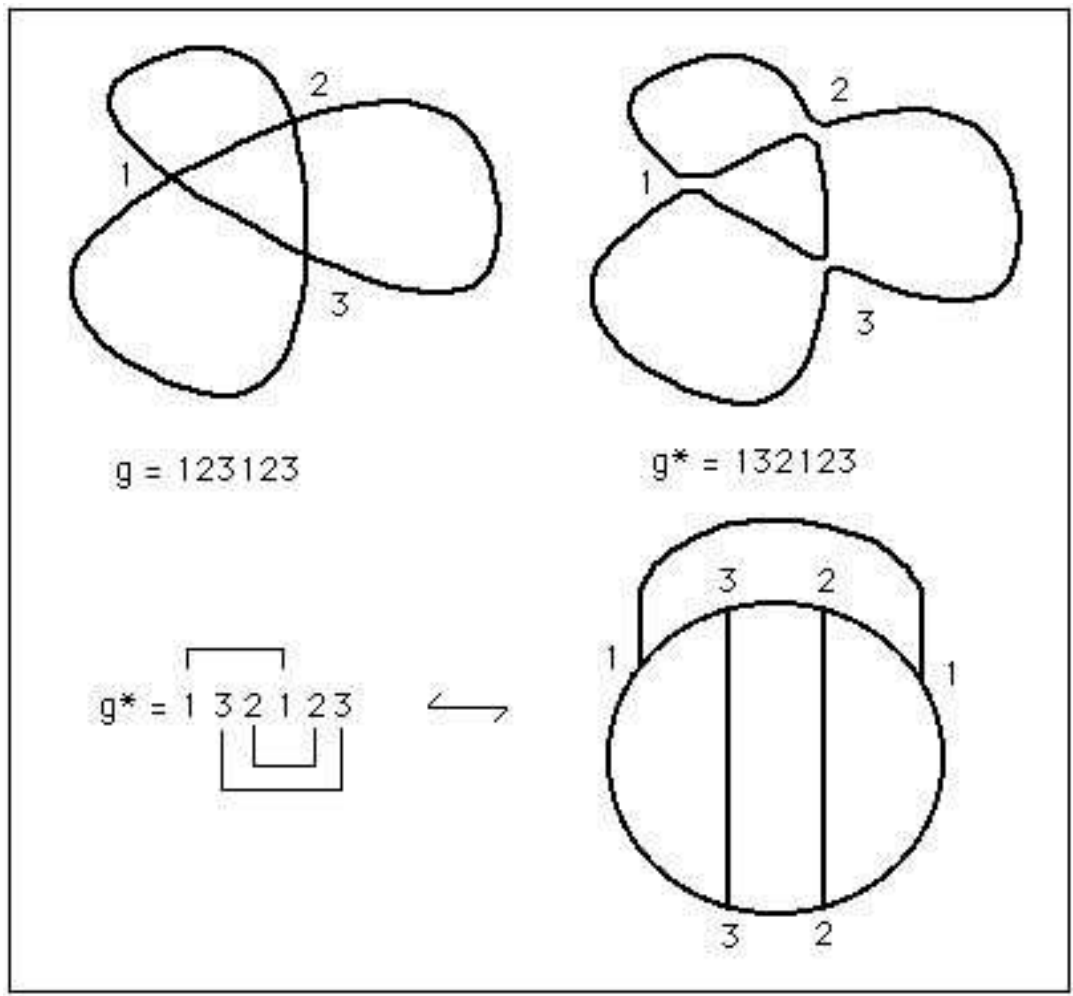
{\bf Figure 7 --- If $g$ is Planar, then $g^{*}$ is Dually Paired.}
\end{center}

\noindent {\bf Remark.} In the case of multiple component codes the algorithm for
constructing $g^{*}$ is modified as follows:  Suppose that the two appearances of
$i$ occur in different components of the code, so that the code to be modified
has the form $$h = i \alpha /i \beta / R$$ where we have  written the two
components as adjacent code segments and started each with $i$ (possible by
rearrangement and cyclic permutation of the segments). Here R denotes the rest of
the code sequences. Then replace h by $$h' = i \alpha i \overline{\beta} / R$$
where $\overline{\beta}$ denotes the rewrite of $\beta$ in reverse order. Note
that  the two components are amalgamated into one as a result of this process. 
Thus after applying this procedure successively to the labels in the code, we
obtain a single code sequence $g^{*}$ from a given multi-component code sequence
$g.$ For example, if $g= 1234/1536/2546$ then we get the following sequence of
partial codes on the way to $g^{*}:$ $$g= 1234/1536/2546$$ $$g' = 12341635/2546 =
23416351/2546$$ $$g'' = 234163512645$$ $$g''' = 236143512645$$ $$g''''=
236146215345$$ $$g'''''= 236146215435$$ $$g^{*} = 236416215435$$ We leave it for
the reader to check that $g^{*}$ is dually paired. \vspace{3mm}

\noindent We have the Lemma below. \vspace{3mm}

\noindent {\bf Lemma 3. } If $g$ is a planar Gauss code, then $g^{*}$ is dually
paired. \vspace{3mm}

\noindent {\bf Proof.} The (easy) proof is omitted. See \cite{R&R}. \vspace{3mm}

Figure 7 illustrates this situation and shows how the desired pairing can be
written directly on the code $g^{*}$ by pairing labels above and below the
typographical line. \vspace{3mm}

\noindent {\bf Lemma 4.} If an evenly intersticed Gauss code $g$ has $g^{*}$
dually paired, then $g$ is the Gauss code of a planar shadow. \vspace{3mm}

\noindent {\bf Proof.} Figure 8 shows how to reconstruct a shadow from any $g$
satisfying the hypotheses of the Lemma. \vspace{3mm}

%\begin{figure}[htbp] \vspace*{90mm} \special{pntg=F8.pntg} \vspace*{13pt}
%\begin{center} { Figure 8 --- Reconstruction from a Gauss Code } \end{center}
%\end{figure} \vspace{3mm}

$$ \picill3inby4in(F8.EPSF) $$
\begin{center}
{\bf  Figure 8 --- Reconstruction from a Gauss Code  }
\end{center}

These lemmas form the essentials of the reconstruction theory for planar Gauss
codes. \vspace{3mm}

\noindent {\bf Definition.}  A Gauss code $g$  is said to be {\em prime} if it
cannot be written as the juxtaposition of two Gauss codes on disjoint collections
of labels.  A non-prime code is said to be {\em composite}. For example, $123123$
is prime but $121234543$ is composite since it is the juxtaposition of $1212$ and
$34543.$ \vspace{3mm}

In reconstructing a shadow from a Gauss code there is a choice of local
orientation of the first crossing in the code.  From then on the local
orientations are determined by the reconstruction algorithm. See Figure 8 for an
example of the procedure.  Once we specify the local orientations in the code,
the corresponding signs of the crossings are determined by whether there is an O
or a U in the code. Thus up to these initial choices of orientation, the signs in
a O/U code are all determined if the code is planar.  It is this result that
gives the proof of Theorem 2. \vspace{3mm}

\noindent {\bf Proof of Theorem 2.}  Note that the reconstruction algorithm will
give a planar embedding for this code with the same local orientations as those
specified in the virtual diagram. In fact, we can assume that the planar
positions of the crossings in the embedded diagram  and the virtual diagram are
identical (up to a global translation if comparison is desired). Now locate those
arcs in the original diagram that involve virtual crossings and move them
one-by-one into the positions indicated by the embedding. To accomplish this,
start at the beginning of the code. Say the code reads $g = a_{1}a_{2} ...
a_{i}a_{i+1} ... a_{n}.$ In the virtual diagram there may be a series of virtual
crossings between $a_{1}$ and $a_{2}$ but there will be no real crossings since
the code is given by $g$. Therefore, the arc from $a_{1}$ to $a_{2}$ can be
replaced (by virtual equivalence) to its position in the embedded diagram.
Continue this process sequentially for $a_{i}a_{i+1}$ and the result is an
equivalence through the virtual category of the original diagram with the
embedded classical diagram. This completes the proof of Theorem 2. //
\vspace{3mm}

\section{Fundamental Group, Crystals, Racks and Quandles}

The fundamental group of the complement of a classical knot can be described by
generators and relations, with one generator for each arc in the diagram and one
relation for each crossing. The relation at a crossing depends upon the type of
the crossing and is either of the form $c=b^{-1} a b$    or    $c=bab^{-1}$  as
shown in Figure 9. \vspace{3mm}

%\begin{figure}[htbp] \vspace*{60mm} \special{pntg=F9.pntg} \vspace*{13pt}
%\begin{center} { Figure 9 --- Generators and Relations for the Fundamental Group}
%\end{center} \end{figure} \vspace{3mm}

$$ \picill3inby2in(F9.EPSF) $$
\begin{center}
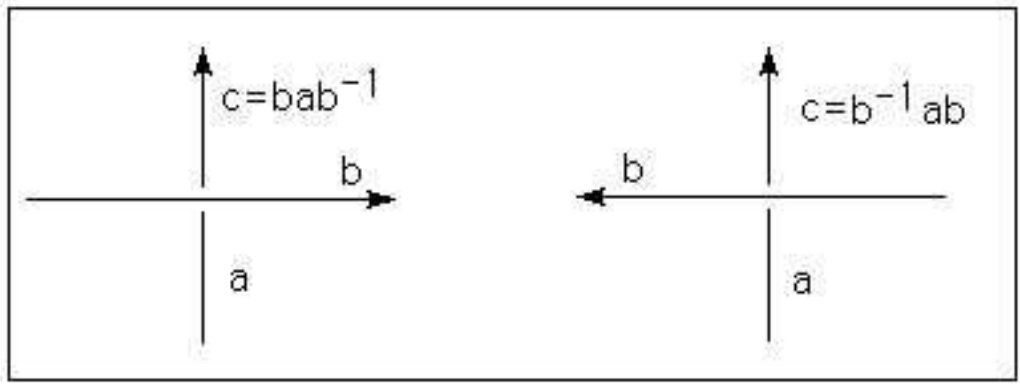
{\bf Figure 9 --- Generators and Relations for the Fundamental Group}
\end{center}

We define the group $G(K)$ of an oriented virtual knot or link by this  same
scheme of  generators and relations. An {\em arc} of a virtual diagram proceeds
from one classical undercrossing to another (possibly the same) classical
undercrossing. Thus no new generators  or relations are added at a virtual
crossing. It is easy to see that $G(K)$ is invariant under all the moves for
virtuals and hence is an invariant of virtual knots. \vspace{3mm}

There are virtual knots that are non-trivial but have trivial fundamental group.
(We say that the fundamental group of a knot is trivial if it is isomorphic to
the infinite cyclic group.)  The virtual $K'$ in Figure 4 is such an example.. We
shall show that $K'$ is a non-trivial virtual in the next section by using a
generalisation of the bracket polynomial. \vspace{3mm}

A generalization of the fundamental group called the quandle, rack or crystal
(depending on notations and history) also assigns relations (in a different
algebra) to each crossing. The quandle generalises to the virtual category.  We
first discuss  the involutory quandle, $IQ(K)$, for a (virtual) knot or link $K.$
The $IQ(K)$ does not depend upon the local orientations of the diagram and it
assigns to each crossing the relation $c=a*b$ as in Figure 10. \vspace{3mm}

%\begin{figure}[htbp] \vspace*{140mm} \special{pntg=F10.pntg} \vspace*{13pt}
%\begin{center} { Figure 10 --- The Involutory Quandle } \end{center} \end{figure}
%\vspace{3mm}

$$ \picill3inby6in(F10.EPSF) $$
\begin{center}
{\bf Figure 10 --- The Involutory Quandle}
\end{center}

The operation $a*b$ is a non-associative binary operation on the underlying set
of the quandle, and it satisfies the following axioms: \vspace{3mm}

\noindent 1.  $a*a=a$ for all a. \vspace{2mm}

\noindent 2.  $(a*b)*b = a$ for all $a$ and $b.$ \vspace{2mm}

\noindent 3. $(a*b)*c = (a*c)*(b*c)$  for all $a$, $b$, $c.$ \vspace{3mm}

The algebra under these axioms with generators and relations as defined above is
called the involutory quandle, $IQ(K).$  It is easy to see that the $IQ(K)$ is
well-defined for $K$ virtual. \vspace{3mm}

An important special case of $IQ(K)$ is the operation  $a*b=2b-a$ where $a$ and
$b$ are elements of a cyclic group $Z/nZ$  for some modulus n. In the case of a
knot $K$, there is a natural choice of modulus $D(K) = Det(M(K))$ where $M(K)$ is
a minor of the matrix of relations associated with the set of equations $c=2b-a.$
 This is called the determinant of the knot, in the classical case, and we shall
call it the determinant of the virtual knot.  If $K$ is virtual then  $|D(K)|$ is
an invariant of $K$.  The virtual knot labelled $K$ in Figure 4 has determinant
equal to 3. The non-triviality of the determinant shows that this knot is knotted
and in fact that it has non-trivial fundamental group. \vspace{3mm}

Another example of an involutory quandle is the operation $a*b = ba^{-1}b.$  In
classical knot theory this yields the fundamental group of the two-fold branched
covering along the knot. \vspace{3mm}

\noindent Here is a useful lemma about the $IQ$ for virtuals. \vspace{3mm}

\noindent {\bf Lemma 5.} $$IQ(K_{vxv})  =  IQ(K_{\overline{x}})$$ where $x$
denotes a crossing in the diagram $K$, $vxv$ denotes that $x$ is flanked by
virtuals, and $K_{\overline{x}}$ denotes the diagram obtained by smoothing the
flanking virtuals, and switching the intermediate crossing. \vspace{2mm}

\noindent In other words the $IQ$ for a classical crossing flanked by two virtual
crossings is the same as the $IQ$ of the diagram where the two virtual crossings
are smoothed and the classical crossing is switched. \vspace{3mm}

\noindent {\bf Proof.}  See Figure 10. // \vspace{3mm}

\noindent {\bf Remark.}  In Figure 10 we illustrate that $IQ(K)=IQ(T)$ where $K$
is the virtual knot also shown in Figure 6 and $T$ is the trefoil knot.
\vspace{3mm}

Finally we discuss the full quandle of a knot and its generalization to virtuals.
 For this discussion the exponential notation of  Fenn and Rourke \cite{F&R} is
convenient.  Instead of $a*b$ we write  $a^{b}$ and assume that there is an
operation of order two

$$a \longrightarrow   \overline{a}$$

\noindent so that

$$\overline {\overline {a}} = a,$$

\noindent and for all $a$ and $b$

$$\overline{a^{b}} = \overline{a}^{b}.$$

\noindent This operation is well-defined for all $a$ in the underlying set $Q$ of
the quandle. \vspace{2mm}

\noindent By definition

$$a^{bc} = (a^{b})^{c}$$ for all $a$,$b$ and $c$ in $Q.$

\noindent The operation of exponentiation satisfies the axioms \vspace{3mm}

\noindent 1. $a^{a} = a$ \vspace{2mm}

\noindent 2. $a^{b \overline{b}} = a$ \vspace{2mm}

\noindent 3. $a^{(b^{c})} = a^{\overline{c} b c}$ \vspace{3mm}

\noindent It follows that the set of the quandle acts on itself by automorphisms

$$ x \longrightarrow x^{a}.$$

\noindent This group of automorphisms is a representation of the fundamental
group of the knot.  Note that if we define $a^{b}$ by the formula

$$a^{b} = bab^{-1}$$

\noindent and

$$\overline{b} = b^{-1},$$

\noindent then we get the fundamental group itself as an example of a quandle. 
The {\em rack} \cite{F&R} or {\em crystal} \cite {K&P} is obtained by eliminating
the first axiom. This makes the rack/crystal an invariant of framed knots and
links. The three axioms correspond to invariance under the three Reidemeister
moves. \vspace{3mm}

If we now compare Lemma 5 with its possible counterpart for the full fundamental
group or the quandle, we see that it no longer holds. Figure 11  shows the new
relations in the quandle that are obtained after smoothing the two virtual
crossings and switching the classical crossing. While the quandle of the
simplified diagram is no longer isomorphic to the original quandle, the fact that
we can articulate the change is often useful in computations. \vspace{3mm}

%\begin{figure}[htbp] \vspace*{120mm} \special{pntg=F11.pntg} \vspace*{13pt}
%\begin{center} { Figure 11 --- Change of Relations for the Full Quandle }
%\end{center} \end{figure} \vspace{3mm}

$$ \picill3inby4.5in(F11.EPSF) $$
\begin{center}
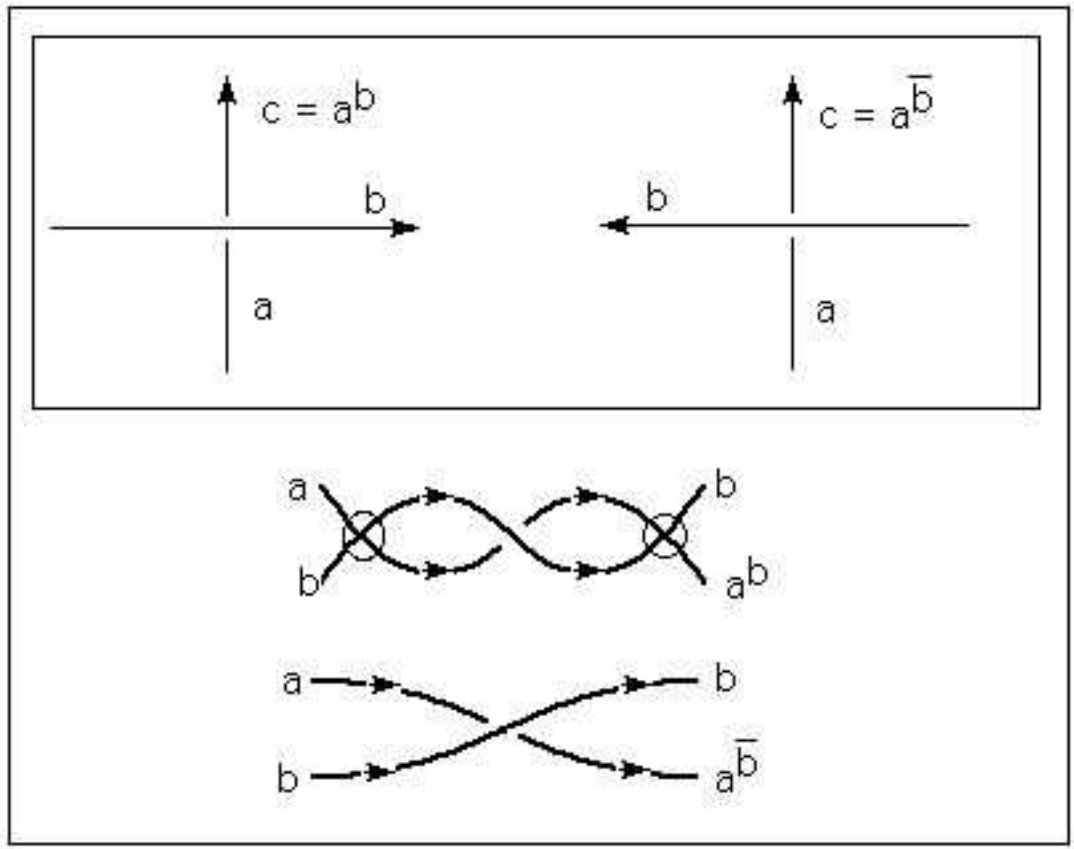
{\bf  Figure 11 --- Change of Relations for the Full Quandle }
\end{center}

\noindent {\bf Example.}  Consider the virtual knot $K$ of Figure 6. We have seen
that $K$ has the same $IQ$ as the trefoil knot. However, the quandle and
fundamental group of $K$ are distinct from those of the trefoil knot, and $K$ is
not equivalent to any classical knot. To see this consider the {\em Alexander
quandle} \cite{K&P} defined by the equations

$$a^{b} = ta +(1-t)b$$

\noindent and

$$a^{\overline{b}} = t^{-1}a +(1-t^{-1})b.$$

\noindent This quandle describes a module (the Alexander module) $M$ over
$Z[t,t^{-1}].$  In the case of the virtual knot $K$ in Figure 6, we have the
generating quandle relations $a^{c} =b$, $b^{a}=c$, $c^{\overline{b}}=a.$ This
results in the Alexander module relations $b = ta + (1-t)c$, $c=tb+(1-t)a$,
$a=t^{-1}c +(1-t^{-1})b.$ From this it is easy to calculate that the module $M(K)
= \{0,m,2m\}$ for a non-zero element $m$ with $3m=0$ and $tm=2m.$ Thus the
Alexander module for $K$ is cyclic of order three. Since no classical knot has a
finite cyclic Alexander module, this proves that $K$ is not isotopic through
virtuals to a classical knot. \vspace{3mm}

Finally, it should be remarked that the full quandle $Q(K)$ classifies a
classical  prime unoriented knot $K$ up to mirror images. By keeping track of a
{\em longitude} for the knot, one gets a complete classification. In the context
of the quandle, the longitude can be described as the automorphism

$$\lambda : Q(K) \longrightarrow Q(K)$$

\noindent defined by the formula

$$\lambda(x) =
x^{a_{1}^{\epsilon_{1}}a_{2}^{\epsilon_{2}}...a_{k}^{\epsilon_{k}}}$$

\noindent where $\{a_{1},a_{2},...,a_{k}\}$ is an ordered list of quandle
generators encountered (as one crosses underneath) as overcrossing arcs as one
takes a trip around the diagram. The $\epsilon$ denotes whether the generator is
encountered with positive or negative orientation, and  $x^{\epsilon}$ denotes
$x$ if $\epsilon = 1$ and $\overline{x}$ if $\epsilon = -1.$  For a given diagram
the longitude is well-defined up to cyclic re-ordering of this list of
encounters.  Exactly the same definition applies to virtual knots. It is no
longer true that the quandle plus longitude classifies a virtual knot, as our
examples of knotted virtuals with trivial fundamental group show. \vspace{3mm}

On the other hand, we can use the quandle to prove the following result. This
proof is due to Goussarov, Polyak and Viro \cite{GPV}. \vspace{3mm}

\noindent {\bf Theorem 6.}  If $K$ and $K'$ are classical knot diagrams such that
$K$ and $K'$ are equivalent under extended virtual Reidemeister moves, then $K$
and $K'$ are equivalent under classical Reidemeister moves. \vspace{3mm}

\noindent {\bf Proof.}  Note that longitudes are preserved under virtual moves
(adding virtual crossings to the diagram does not change the expression for a
longitude).  Thus an isomorphism from $Q(K)$ to $Q(K')$ induced by extended moves
preserves longitudes.  Since the isomorphism class of the quandle plus longitudes
classifies classical knots, we conclude that $K$ and $K'$ are classically
equivalent. This completes the proof. // \vspace{3mm}

\noindent {\bf Remark.} We would like to see a purely combinatorial proof of
Theorem 6.

\subsection{The GPV Example and a Generalisation.}

We end this section with a variation of an example \cite{Viro}, \cite{GPV} that
shows that it is possible to have a virtual knot $K$ with $Q(K)$ not isomorphic
with $Q(K^{*})$ where $K^{*}$ is the mirror image of $K.$ In other words, there
are {\em two} quandles,  or two fundamental  groups associated with any given
virtual knot! \vspace{3mm}

This example is a slightly different take on an observation in \cite{GPV}. Let
$K$ be a given (virtual) diagram, drawn in the plane. Pick the diagram up and
turn it over (note that the crossings change diagrammatically, but correspond to
the result of physically turning over the layout of criss-crossing strands with
welds at the virtual crossings). Let $Flip(K)$ denote this overturned diagram.
Define a new quandle $Q^{*}(K)$ by the formula $Q^{*}(K) = Q(Flip(K)).$ 
Goussarov, Polyak and Viro take their "other" fundamental group to be the one
defined by generators and relations obtained by "looking at the knot from the
other side of the plane."  At the quandle level this is the same as taking
$Q^{*}(K).$ It is easy to see that {\em $Q^{*}(K)$ is isomorphic to $Q(K^{*}).$}
(Just note that if $c=a^{\overline{b}}$ then $\overline{c} =
\overline{a^{\overline{b}}} = \overline{a}^{\overline{b}}.$ Use this to check
that the two quandles are isomorphic through  the mapping $a \longrightarrow
\overline{a}$ taking one to the other.) Thus our version of this example is
mathematically equivalent to the GPV version. \vspace{3mm}

%\begin{figure}[htbp] \vspace*{120mm} \special{pntg=F12.pntg} \vspace*{13pt}
%\begin{center} { Figure 12 ---  $Q(K)$ is distinct from $Q(K^{*}).$ }
%\end{center} \end{figure} \vspace{3mm}

$$ \picill3inby5in(F12.EPSF) $$
\begin{center}
{\bf Figure 12 ---  $Q(K)$ is distinct from $Q(K^{*}).$ }
\end{center}

In Figure 12 the reader will find $K$ and $K^{*}$ with labelled arcs
$a$,$b$,$c$,$d.$ In $K$ the quandle relations are  $a = b^{d}$, $b=c^{d}$,
$c=d^{b}$, $d=a^{b}.$  The three-coloring of $K$ in $Z/3Z$  with $a=0$, $b=2$,
$c=0$, $d=1$ demonstrates that this quandle, and hence the fundamental group of
$K$ is non-trivial. On the other hand, $K^{*}$ has quandle relations $a=a^{a}$,
$c=b^{c}$, $d=c^{c}$, $a=d^{a}$, giving a trivial quandle. Thus $K$ is
distinguished from $K^{*}$ by the quandle.  This example also shows that $K$ has
non-trivial Alexander polynomial (using the fundamental group to define the
Alexander polynomial - there is more than one Alexander polynomial for virtuals)
but $K^{*}$ has Alexander polynomial equal to $1.$ \vspace{3mm}

%\begin{figure}[htbp] \vspace*{60mm} \special{pntg=F13.pntg} \vspace*{13pt}
%\begin{center} { Figure 13 --- The 1--1 Tangle $W$ } \end{center} \end{figure}
%\vspace{3mm}

$$ \picill3inby4in(F13.EPSF) $$
\begin{center}
{\bf Figure 13 ---- The 1--1 Tangle $W$ }
\end{center}

We generalise this example by considering the 1--1 tangle $W$ shown in Figure 13.
Replacing a straight arc in a knot diagram by $W$ does not affect the quandle,
while replacing by its mirror image $W^{*}$ changes the quandle relations in a
generally non-trivial way (as in the GPV example). Insertion of $W$ into knot
diagrams produces infinitely many examples of pairs of virtual knots with the
same quandle but different Jones polynomials. This last statement will be
verified in the next section. \vspace{3mm}

\section{Bracket Polynomial and Jones Polynomial}

The bracket polynomial  \cite{KState} extends to virtual knots and links by
relying on the usual formula for the state sum of the bracket, but allowing the
closed loops in the state to have virtual intersections. Each loop is still
valued at $d=-A^{2} -A^{-2}$  and the expansion formula \vspace{3mm}

$$<K>  = A<K_{a} > + A^{-1}<K_{b}>$$

\noindent still holds where $K_{a}$ and $K_{b}$ denote the result of replacing a
single crossing in $K$ by smoothings of type $a$ and type $b$ as illustrated in
Figure 14. \vspace{3mm}

%\begin{figure}[htbp] \vspace*{70mm} \special{pntg=F14.pntg} \vspace*{13pt}
%\begin{center} { Figure 14 --- Bracket Smoothings } \end{center} \end{figure}
%\vspace{3mm}

$$ \picill3inby2in(F14.EPSF) $$
\begin{center}
{\bf Figure 14 --- Bracket Smoothings }
\end{center}

We must check that this version of the bracket polynomial is invariant under all
but the  first Reidemeister move (See the moves shown  in Figure 2 ).  Certainly
the usual arguments apply to the moves of type (A).  Moves of type (B) do not
disturb the loop counts and so leave bracket invariant. Finally the move of type
(C) receives the verification illustrated in Figure 15. This completes the proof
of the invariance of the generalised bracket polynomial under move (C).
\vspace{3mm}

%\begin{figure}[htbp] \vspace*{140mm} \special{pntg=F15.pntg} \vspace*{13pt}
%\begin{center} { Figure 15 --- Type (C) Invariance of the Bracket } \end{center}
%\end{figure} \vspace{3mm}

$$ \picill3inby5in(F15.EPSF) $$
\begin{center}
{\bf Figure 15 --- Type (C) Invariance of the Bracket }
\end{center}

We define the writhe $w(K)$ for an oriented virtual to be the sum of the crossing
signs---just as in the classical case. \vspace{3mm}

The $f$-polynomial is defined by the formula

$$f_{K}(A) = (-A^{3})^{-w(K)} <K>(A).$$

\noindent The Laurent polynomial, $f_{K}(A)$ is invariant under all the virtual
moves including the classical move of type I. \vspace{3mm}

\noindent {\bf Remark.} It is worth noting that $f_{K}$ can be given a state
summation of its own. Here we modify the vertex weights of the bracket state sum
to include a factor of $-A^{-3}$ for each crossing of positive sign, and a factor
of $A^{+3}$ for each factor of negative sign.  It is then easy to see that

$$f_{K_{+}} = -A^{-2} f_{K_{0}} - A^{-4} f_{K_{\infty}}$$ $$f_{K_{-}} = -A^{+2}
f_{K_{0}} -A^{+4}  f_{K_{\infty}}$$

\noindent where $K_{+}$ denotes $K$ with a selected positive crossing, $K_{-}$
denotes the result of switching only this crossing, $K_{0}$ denotes the result of
making the oriented smoothing of this crossing, and $K_{\infty}$  denotes the
result of making an unoriented smoothing at this crossing.  The states in this
oriented state sum acquire  sites with unoriented smoothings, but the procedure
for evaluation is the same as before. For each state  we take the product of the
vertex weights multiplied by  $d^{||S||-1}$  where $d=-A^{2} - A^{-2}$  and
$||S||$ denotes the number of loops in the state. Then $f_{K}$ is the sum of
these products, one for each state. \vspace{3mm}

The following Lemma makes virtual calculations easier. \vspace{3mm}

\noindent {\bf Lemma 7.}   $< K_{vxv}>   =  < K_{x} >$  where $x$ denotes a
crossing in the diagram $K$, $vxv$ denotes that $x$ is flanked by virtuals, and
$K_{x}$ denotes the diagram obtained by smoothing the flanking virtuals, and
leaving the crossing the same. \vspace{3mm}

\noindent {\bf Proof.}    The proof is shown in Figure 16.//

%\begin{figure}[htbp] \vspace*{120mm} \special{pntg=F16.pntg} \vspace*{13pt}
%\begin{center} { Figure 16 --- Removal of Flanking Virtual Crossings }
%\end{center} \end{figure} \vspace{3mm}

$$ \picill3inby3.5in(F16.EPSF) $$
\begin{center}
{\bf Figure 16 --- Removal of Flanking Virtual Crossings  }
\end{center}

Note that this result has the opposite form from our corresponding Lemma about
the involutory quandle $IQ(K).$   As a result we get an example of a virtual knot
that is non-trivial (via the $IQ$) but has $f_{K} =1.$  {\em Hence we have a
virtual knot $K$ with Jones polynomial equal to 1.}   The example is shown in
Figure 17. Note that in Figure 10 we illustrated that  this  $K$ has the same
involutory quandle as the trefoil knot. We will see in  Section 6 that $K$ is not
equivalent to a classical knot. \vspace{3mm}

%\begin{figure}[htbp] \vspace*{60mm} \special{pntg=F17.pntg} \vspace*{13pt}
%\begin{center} { Figure 17 --- A Knotted Virtual  with Trivial Jones Polynomial }
%\end{center} \end{figure} \vspace{3mm}

$$ \picill3inby2in(F17.EPSF) $$
\begin{center}
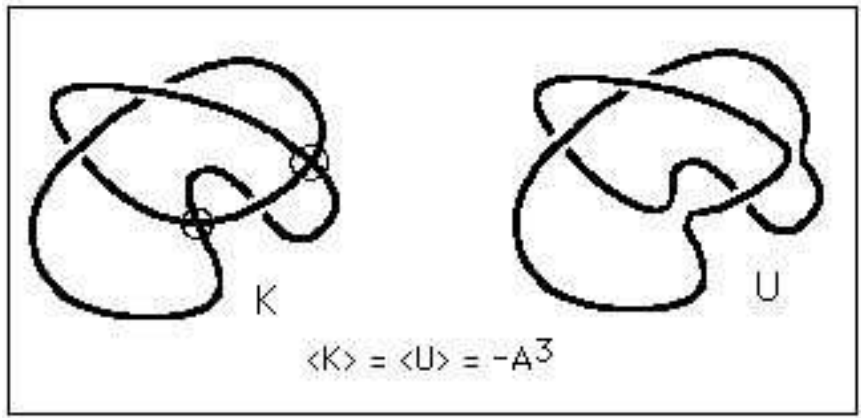
{\bf Figure 17 --- A Knotted Virtual  with Trivial Jones Polynomial  }
\end{center}

We now compute the bracket polynomial for our previous example with trivial
fundamental group and we find that $<K'> = A^{2} + 1 - A^{-4}$  and $f_{K'} =
(-A^{3})^{-2}<K'> = A^{-4} + A^{-6} - A^{-10}.$  Thus $K'$ has a non-trivial
Jones polynomial. See Figure 18. 
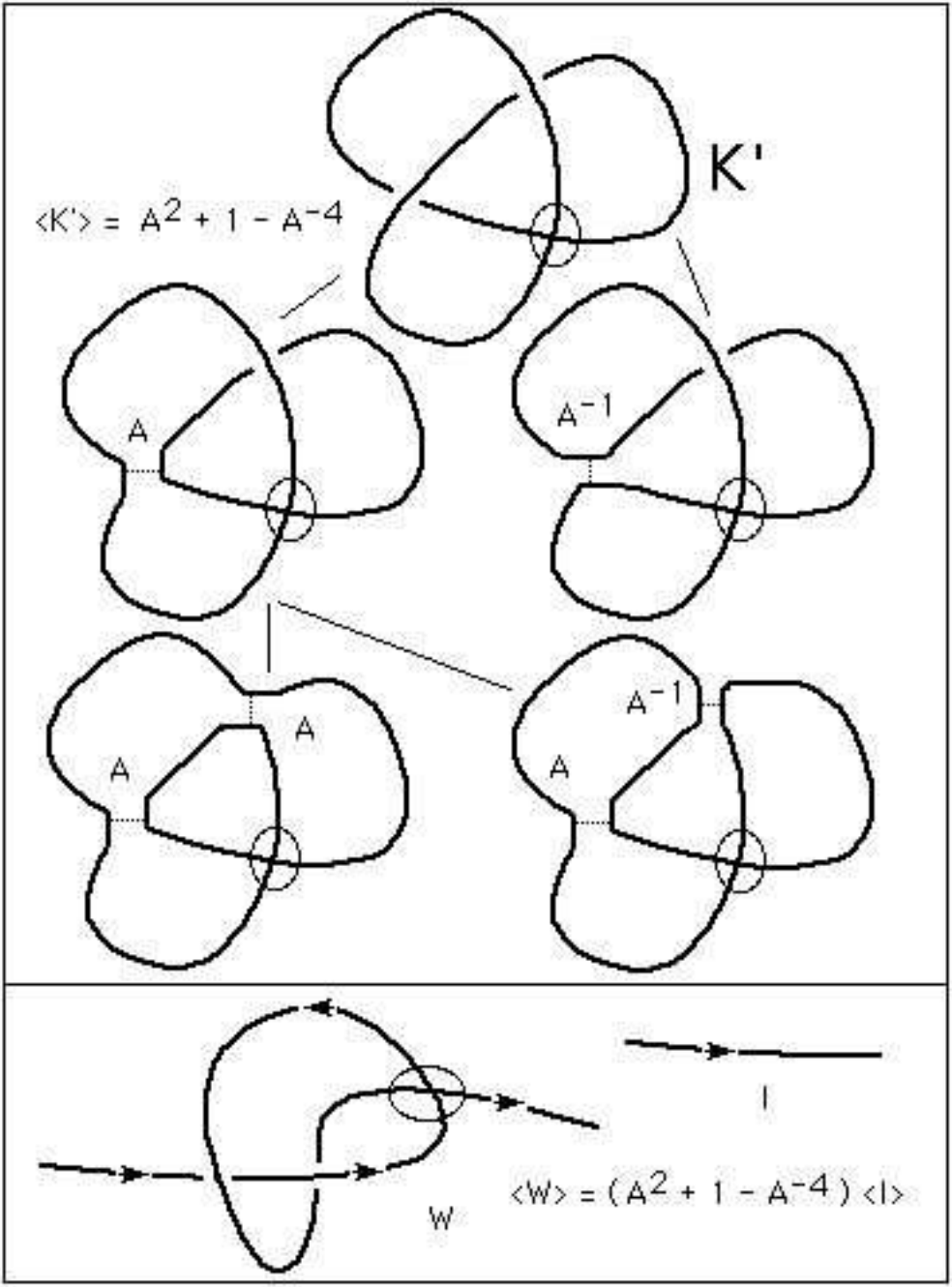
\vspace{3mm}

%\begin{figure}[htbp] \vspace*{140mm} \special{pntg=F18.pntg} \vspace*{13pt}
%\begin{center} { Figure 18 --- Calculation of $<K'>$ } \end{center} \end{figure}
%\vspace{3mm}

$$ \picill3inby7in(F18.EPSF) $$
\begin{center}
{\bf Figure 18 --- Calculation of $<K'>$ }
\end{center}

In Figure 18 we also indicate the result of placing the tangle W, discussed in
Figure 13,  into another knot or link.  Since this is the same as taking a
connected sum with $K'$ it has the effect of multiplying the bracket polynomial
by $A^{2} + 1 - A^{-4}.$  Thus if $L$ is any knot or link and $K' + L$ denotes
the connected sum of $K'$ along some component of L then $<K' + L> =  (A^{2} + 1
- A^{-4}) <L>$  while $Q(K' + L) = Q(L)$  (as we verified in the last section). 
Thus for any knot L, successive connected sums with $K'$ produces an infinite
family of distinct virtual knots, all having the same quandle (hence same
fundamental group). \vspace{3mm}

Finally, we note that if the knot is given as embedded in $S_{g} \times I$ for a
surface of genus g, and if its virtual knot diagram $K$ is obtained by projecting
the diagram on $S_{g}$ into the plane, then $<K>$ computes the value of the
extension of the bracket to the knots in $S_{g} \times I$  where all the loops
have the same value $d = -A^{2} -A^{-2}.$  This is the first order bracket for
link diagrams on a surface. \vspace{3mm}

In Figure 4 we illustrated the non-trivial knot K with trivial Jones polynomial
as embedded in $S_{1} \times I.$  This knot in $S_{1} \times I$ is actually not
trivial as can be seen from the higher Jones polynomials that discriminate loops
in different isotopy classes on the surface. \vspace{3mm}

%\begin{figure}[htbp] \vspace*{120mm} \special{pntg=F19.pntg} \vspace*{13pt}
%\begin{center} { Figure 19 --- A Knot in $S_{1} \times I$ with Trivial Jones
%Polynomial} \end{center} \end{figure} \vspace{3mm}

$$ \picill3inby4in(F19.EPSF) $$
\begin{center}
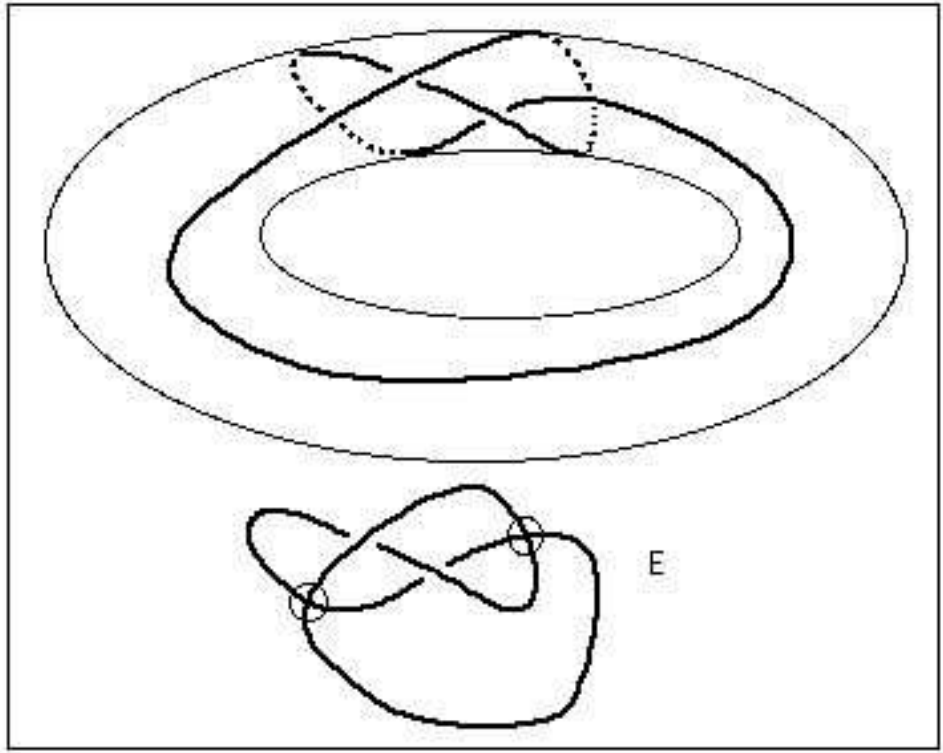
{\bf Figure 19 --- A Knot in $S_{1} \times I$ with Trivial Jones
Polynomial}
\end{center}

In Figure 19  is another example of a virtual knot $E$ and a corresponding
embedding in $S_{1} \times I$.  In this case, $E$ is a trivial virtual knot (as
is shown in Figure 3), but the embedding of $E$ in $S_{1} \times I$ is
non-trivial (even though it has trivial fundamental group and trivial bracket
polynomial).  The non-triviality of this embedding is seen by simply observing
that it carries a non-trivial first homology class in the thickened torus. In
fact, if you expand the state sum for the bracket polynomial and keep track of
the isotopy classes of the curves in the states, then the bracket calculation
also shows this non-triviality by exhibiting as its value a single state with a
non-contractible curve. \vspace{3mm}

Virtual knot theory provides a convenient calculus for working with knots in
$S_{g} \times I.$ The virtuals carry many properties of knots in $S_{g} \times I$
that are independent of the choice of embedding and genus.  This completes our
quick survey of the properties of the bracket polynomial and Jones polynomial for
virtual knots and links. Just as uncolorable graphs appear when one goes beyond
the plane (for planar graph coloring problems) so knots of unit Jones polynomial
appear as we leave the diagrammatic plane into the realm of the Gauss codes.
\vspace{3mm}

\section{Quantum Link Invariants}

There are virtual link invariants corresponding to every quantum link invariant
of classical links. However this must be said with a caveat: We do not assume
invariance under the first classical Reidemeister move (hence these are
invariants of regular isotopy) and we do not assume invariance under the flat
version of the first Reidemeister move in the (B) list of virtual moves. 
Otherwise the usual tensor or state sum formulas for quantum link invariants
extend to this generalized notion of regular isotopy invariants of virtual knots
and links.  In this section we illustrate this method by taking a different
generalisation of the bracket that includes virtual framing. We apply this new
invariant to distinguish a virtual knot that has Jones polynomial equal to one
and a trivial fundamental group. \vspace{3mm}

In order to carry out this program, we quickly recall how to construct quantum
link invariants in the unoriented case. See \cite{K&P} for more details. The link
diagram is arranged with respect to a given "vertical" direction in the plane so
that perpendicular lines to this direction intersect the diagram transversely or
tangentially at maxima and minima.. In this way the diagram can be seen as
constructed from a pattern of interconnected maxima, minima and crossings---as
illustrated in Figure 20. 
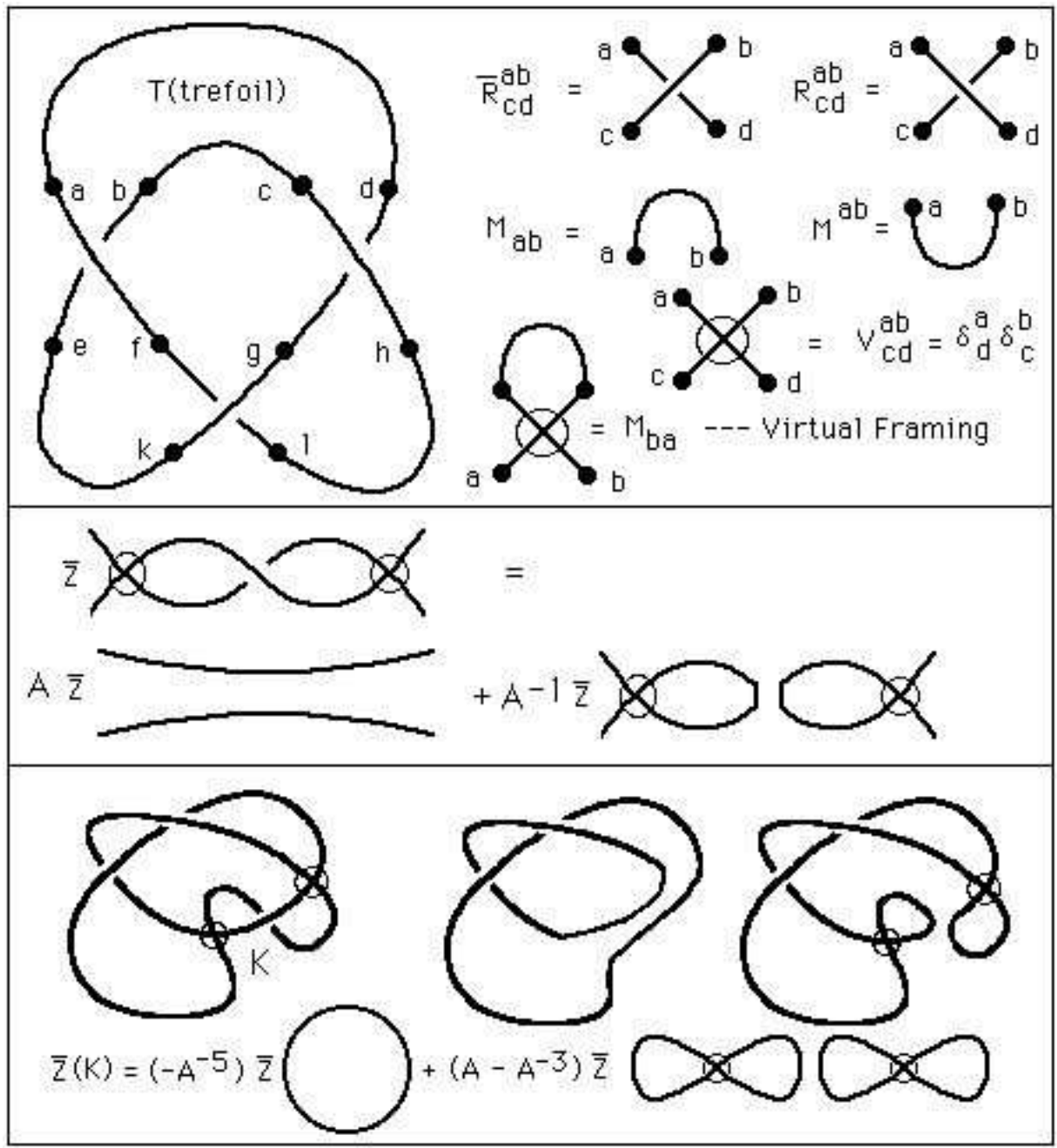
\vspace{3mm}

%\begin{figure}[htbp] \vspace*{160mm} \special{pntg=F20.pntg} \vspace*{13pt}
%\begin{center} { Figure 20 --- Quantum Link Invariants } \end{center}
%\end{figure} \vspace{3mm}

$$ \picill3inby7in(F20.EPSF) $$
\begin{center}
{\bf Figure 20 --- Quantum Link Invariants}
\end{center}

As illustrated in Figure 20, we associate symbols $M_{ab}$ and $M^{ab}$ to minima
and maxima respectively, and symbols $R^{ab}_{cd}$ and $\overline{R}^{ab}_{cd}$
to the two types of crossings. The indices on these symbols indicate how they are
interconnected. Each maximum or minimum has two lines available for connection
corresponding to the indices $a$ and $b.$ Each $R$ , $\overline{R}$  has four
lines available for connection.  Thus the symbol sequence

$$T(K) = M_{ad}M_{bc}M^{ek}M^{lh}R^{ab}_{ef}R^{cd}_{gh}\overline{R}^{fg}_{kl}$$

\noindent represents the trefoil knot as shown in Figure 20.  Since repeated
indices show the places of connection, there is no necessary order for this
sequence of symbols. I call $T(K)$ an {\em abstract tensor} expression for the
trefoil knot $K$. \vspace{3mm}

By taking matrices (with entries in a commutative ring) for the $M$'s and the
$R$'s it is possible to re-interpret the abstract tensor expression as a
summation of products of matrix entries over all possible choices of indices in
the expression.  Appropriate choices of matrices give rise to link invariants. 
If $K$ is a knot or link and $T(K)$ its associated tensor expression, let $Z(K)$
denote the evaluation of the tensor expression that corresponds to the above
choice of matrices. We will assume that the matrices have been chosen so that
$Z(K)$ is an invariant of regular isotopy. \vspace{3mm}

The generalisation of the quantum link invariant $Z(K)$ to virtual knots and
links is quite straightforward. We simply ignore the virtual crossings in the
diagram. Another way to put this is that we take each virtual crossing to be
represented by crossed Kronecker deltas as in Figure 20.  The virtual crossing is
represented by the tensor  $$V^{ab}_{cd} = \delta^{a}_{d} \delta^{b}_{c}.$$ Here
$\delta^{a}_{b}$ is the Kronecker delta. It is equal to $1$ if $a=b$ and is equal
to $0$ otherwise. (Note that the Kronecker delta is well defined as an abstract
tensor.) \vspace{3mm}

In extending $Z(K)$ to virtual knots and links  by this method we cannot hope to
obtain invariance under the type I virtual move. In fact, as Figure 20 shows, the
presence of a virtual curl is indexed by the transpose $M_{ba}$ of the tensor
$M_{ab}$. Thus we define {\em virtual regular isotopy} to be invariance under all
the extended Reidemeister moves for virtuals except type (A)I and (B)I.  It is
easy to see that $Z(K)$ extends in this way when $Z(K)$ is an invariant of
regular isotopy for classical links. \vspace{3mm}

In particular the bracket polynomial for classical knots is obtained by letting
the indices run over the set $\{1,2\}$ with $M^{ab} = M_{ab}$ for all $a$ and $b$
and $M_{11}=M_{22} = 0$ while $M_{12} = iA$ and $M_{21} = -iA^{-1}$ where $i^{2}
= -1.$  The R's are defined by the equations

$$R^{ab}_{cd} = A M^{ab}M_{cd} + A^{-1} \delta^{a}_{c} \delta^{b}_{d}$$

$$\overline{R}^{ab}_{cd} = A^{-1} M^{ab}M_{cd} + A \delta^{a}_{c}
\delta^{b}_{d}$$

\noindent These equations for the $R$'s are the algebraic translation of the
smoothing identities for the bracket polynomial.    Then we have \vspace{3mm}

\noindent {\bf Theorem 8.}  With $Z(K)$ defined as above and $K$ a classical knot
or link, $Z(K) = d<K>$ where $d = -A^{2} - A^{-2}$. \vspace{3mm}

\noindent {\bf Proof.} See \cite{K&P}.// \vspace{3mm}

For this extension of $Z(K)$  to virtuals there is a state summation similar to
that of the bracket polynomial. For this let $C$ be a diagram in the plane that
has only virtual crossings. View this diagram as an immersion of a circle in the
plane. Let $rot(C)$ denote the absolute value of the Whitney degree of $C$ as a
immersion in the plane. (Since $C$ is unoriented only the absolute value of the
Whitney degree is well-defined.). The Whitney degree of an oriented plane
immersion is the total algebraic number of $2 \pi$ turns of the unit tangent
vector to the curve as the curve is traversed once. Let $d(C)$ be defined by the
equation $$d(C) = (-1)^{rot(C)} (A^{2 rot(C)} + A^{-2 rot(C)}).$$ Let $S$ be a
state of a virtual diagram $K$ obtained by smoothing each classical crossing in
$K.$  Let $C \in S$ mean that $C$ is one of the curves in $S.$  Let $<K|S>$
denote the usual product of vertex weights ($A$ or $A^{-1}$) in the bracket state
sum.  Then \vspace{3mm}

\noindent {\bf Proposition 9.}  The invariant of virtual regular isotopy $Z(K)$
is described by the following state summation.

$$Z(K) = \sum_{S} <K|S> \prod_{C \in S} d(C)$$

\noindent where the terms in this formula are as defined above. Note that $Z(K)$
reduces to $d<K>$ when $K$ is a classical diagram. \vspace{3mm}

\noindent {\bf Proof.} The proof is a calculation based on the tensor model
explained in this section. The details of this calculation are omitted.//
\vspace{3mm}

\noindent {\bf Remark.}  The state sum in Proposition 9 generalises to an
invariant of virtual regular isotopy with an infinite number of polynomial
variables, one for each regular homotopy class of unoriented curve immersed in
the plane. To make this generalisation, let $A_{n}$ for $n=0,1,2,3,...$ denote a
denumerable set of commuting independent variables.  If $C$ is an immersed curve
in the plane, define $Var(C) = A_{n}$ where $n = rot(C)$, the absolute value of
the Whitney degree of $C.$  We take $A_{1} = -A^{2} - A^{-2}$ as before, but the
other variables are independent of each other and  of $A$. \vspace{3mm}

\noindent Now define the generalisation of $Z(K)$, denoted $\overline{Z}(K),$ by
the formula

$$\overline{Z}(K) = \sum_{S} <K|S> \prod_{C \in S} Var(C).$$

In this definition we have replaced the evaluation $d(C)$ by the corresponding
variable $Var(C).$  In Figure 20 we illustrate the result of calculating
$\overline{Z}(K)$  for a knot $K$ with unit Jones polynomial. The result is

$$\overline{Z}(K) = (-A^{-5})A_{1} + (A-A^{-3})A_{0}^{2}.$$

Since the coefficients of $A_{1}$ and $A_{0}$  are themselves invariants of
virtual regular isotopy it follows, as we already knew, that $K$ is a non-trivial
virtual. This non-triviality is detected by our refinement $\overline{Z}(K)$ of
the bracket polynomial. A similar phenomenon of refinement of invariants happens
with other quantum link invariants.  This will be the subject of a separate
paper. \vspace{3mm}

We end this section with an application of Proposition 9. Let $D$ be the virtual
knot diagram shown in Figure 21. \vspace{3mm}

%\begin{figure}[htbp] \vspace*{140mm} \special{pntg=F21.pntg} \vspace*{13pt}
%\begin{center} { Figure 21 --- $D$ has Unit Jones Polynomial and Trivial
%Fundamental Group } \end{center} \end{figure} \vspace{3mm}

$$ \picill3inby6in(F21.EPSF) $$
\begin{center}
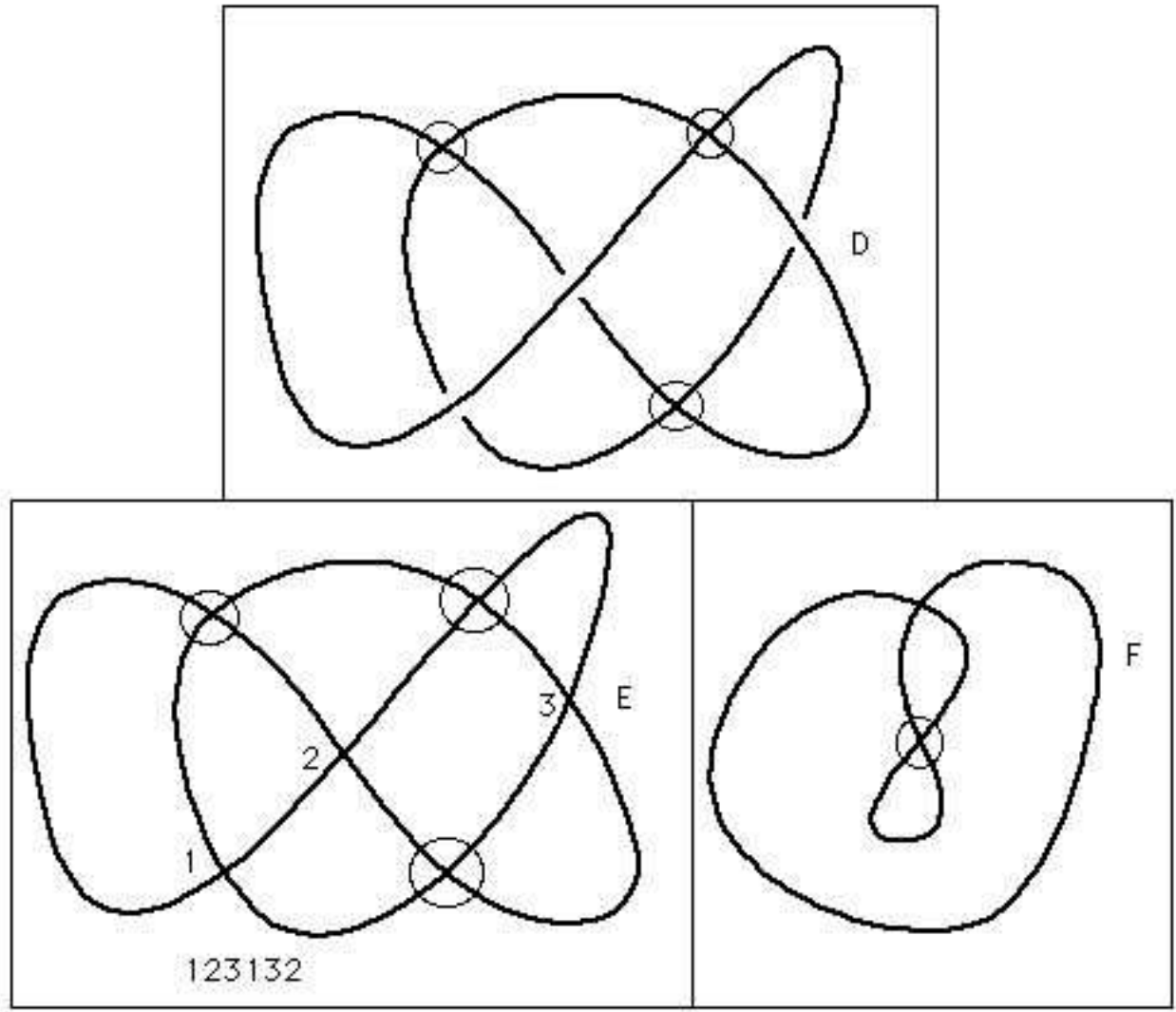
{\bf Figure 21 --- $D$ has Unit Jones Polynomial and Trivial
Fundamental Group }
\end{center}

\noindent It is easy to see that $D$ has $f$-polynomial equal to $1$, hence Jones
polynomial equal to $1.$ Use Lemma 7 to show that $<D> = -A^{3}.$  $D$ also has
trivial fundamental group and quandle. Is $D$ a non-trivial virtual knot? The
answer is yes! It follows from the calculation of $Z(D).$ We omit the
calculation, but give the result

$$Z(D) = A^{7} -A^{5}-4A^{3}+2A+A^{-1}-A^{-3}.$$

\noindent It follows from this that $D$ cannot be regularly isotopic to a
standard virtual curl form. Hence $D$ must be virtually knotted in the regular
isotopy category.  On the other hand, I do not yet have a proof that $D$ is
virtually knotted under the original definition that allows the addition and
removal of virtual framings. This example shows both the power and limitation of
using the quantum invariants to study virtual knots. \vspace{3mm}

There should be a direct way to see that $D$ is knotted. Let $E$ denote the {\em
shadow} of the diagram $D.$ That is, replace the classical crossings in $D$ with
flat crossings. Regard the flat crossings as {\em distinct} from virtual
crossings, so that we get the rules for virtual isotopy of flat diagrams shown in
Figure 22. By these rules a flat diagram corresponds to an oriented Gauss code
without over or under crossing specifications. The virtual moves preserve the
Gauss code just as before. \vspace{3mm}

%\begin{figure}[htbp] \vspace*{140mm} \special{pntg=F22.pntg} \vspace*{13pt}
%\begin{center} { Figure 22 --- Flat Virtual Moves} \end{center} \end{figure}
%\vspace{3mm}

$$ \picill3inby6in(F22.EPSF) $$
\begin{center}
{\bf Figure 22 --- Flat Virtual Moves}
\end{center}

\noindent $E$ is illustrated in Figure 21. Is $E$ flat virtually equivalent to a
circle with curls and virtual curls?  I conjecture that the answer is no. Simpler
examples of this sort of irreducibility are easy to produce. The diagram $F$ in
Figure 21 is irreducible because $\overline{Z}(F)(1) = -2A_{0} + A_{2}$, as is
easy to compute. \vspace{3mm}

These examples lead us to the following definition. We call the {\em shadow code}
of a diagram the underlying Gauss code of that diagram without any specifications
of orientation or over/under crossing.  We say that a virtual diagram is {\em
almost classical} if its shadow code is planar. \vspace{3mm}

\noindent {\bf Conjecture.} There does not exist a non-trivial  almost classical
virtual knot with both trivial fundamental group and trivial Jones polynomial.

\section{Virtual Vassiliev Invariants}

We now study embeddings into $R^{3}$  (Euclidean three-space) of 4-valent  graphs
up to {\em rigid vertex isotopy.}  In rigid vertex isotopy one can think of each
graphical vertex as a rigid disk. The four graphical edges incident to this
vertex are attached to the boundary of the disk at four specific points.  In a
rigid vertex isotopy the embedded edges of the graph can be isotoped freely, but
the disks must move without deformation in the course of the isotopy.  A
consequence of this definition \cite{KaufDiag} is that the diagrammatic moves
shown in Figure 23 capture rigid vertex isotopy just as the Reidemeister moves
capture ambient isotopy. Figure 23 shows only the move types that are added to
the usual list of Reidemeister moves.  Two graphs with diagrammatic projections
$G_{1}$ and $G_{2}$ are rigid vertex isotopic if and only if there is a series of
moves of this type  joining the two diagrams. \vspace{3mm}

%\begin{figure}[htbp] \vspace*{100mm} \special{pntg=F23.pntg} \vspace*{13pt}
%\begin{center} { Figure 23 --- Moves for Rigid Vertex Embeddings} \end{center}
%\end{figure} \vspace{3mm}

$$ \picill3inby4in(F23.EPSF) $$
\begin{center}
{\bf Figure 23 --- Moves for Rigid Vertex Embeddings}
\end{center}

\noindent In Figure 23 there is also illustrated the one addition to virtual
moves that is needed to complete the move set for rigid vertex isotopy of virtual
knotted graphs. In this addition an arc with consecutive virtual crossings is
moved to a new position across a rigid vertex.  Here we must make a distinction
between the graphical rigid vertices and the virtual vertices in the diagrams. 
Once this is done we directly extend discussions of invariants of rigid vertex
graphs to invariants of virtual rigid vertex graphs.  To see how this is done we
will discuss invariants obtained by insertion into the vertices of a graph. We
shall always mean the extension to virtual equivalence when we refer to ambient
isotopy or to rigid vertex isotopy. \vspace{3mm}

If one replaces each node of a (virtual) rigid vertex graph $G$ with a tangle
(possibly virtual) to form a virtual link $K$ , then any rigid equivalence of $G$
induces a corresponding equivalence of $K$.  The consequence of this remark is
that we can obtain invariants of rigid vertex graphs from any invariant of
(virtual) knots and links by taking a systematic choice of tangle insertion. 
That is, if we have chosen tangle insertions $T_{1}$, ..., $T_{n}$ let $\{v_{1},
..., v_{m} \}$ denote the set of vertices of $G$ and let $a=(a_{1}, ..., a_{m})$
with  $1 \leq a_{i} \leq n$ denote a choice of tangle insertion for each vertex
of $G.$  Then let $G^{a}$ denote the result of inserting tangle $T_{a_{i}}$ at
node $i$ in $G.$  Suppose that $R(K)$ is an ambient isotopy invariant of virtual
knots and links $K.$  Then define an extension of $R$ to graphical imbeddings by
the formula

$$R(G) = \sum_{a} x_{a_{1}} ... x_{a_{m}} R(G^{a})$$

\noindent where $\{   x_{j} | j = 1 ,..., m \}$  is a new set of variables (or
constants) independent of the variables already present in the invariant $R.$ Our
discussion shows that $R(G)$ is an invariant of virtual graph embeddings $G.$
\vspace{3mm}

While it is of interest to explore this larger class of induced invariants, we
shall restrict ourselves to the generalisation of {\em Vassiliev invariants}.  A
Vassiliev invariant $v$  is an invariant of rigid vertex (virtual) 4-valent
graphs  that satisfies

$$v(G|*) = v(G|+) - v(G|-)$$

\noindent where $(G|*)$  denotes an oriented graph $G$ with a chosen vertex $*.$ 
$(G|+)$ denotes the result of replacing the vertex $*$ with a positive crossing
and $(G|-)$ is the result of replacing it with a negative crossing. See Figure
24. This is the traditional definition of a Vassiliev invariant and we adopt it
verbatim for virtuals. \vspace{3mm}

%\begin{figure}[htbp] \vspace*{60mm} \special{pntg=F24.pntg} \vspace*{13pt}
%\begin{center} { Figure 24 --- Vassiliev Invariant Identity } \end{center}
%\end{figure} \vspace{3mm}

$$ \picill3inby2.5in(F24.EPSF) $$
\begin{center}
{\bf Figure 24 --- Vassiliev Invariant Identity}
\end{center}

\noindent {\bf Definition.}  Let $N(G)$ denote the number of vertices in the
4-valent graph $G.$  A Vassiliev invariant $v$  is said to be of {\em graphical
finite type n} if $v(G) = 0$ whenever $N(G) > n.$ Note that this definition says
nothing about the number of virtual crossings in the graph $G.$ \vspace{3mm}

Useful examples of virtual invariants of graphical finite type are obtained by
taking the coefficients of $x^{m}$ in

$$F_{K}(x) = f_{K}(e^{x})$$

\noindent where $F_{K}(x)$ is extended to 4-valent graphs by the difference
formula

$$F_{K}(x)(G|*) = F_{K}(x)(G|+) - F_{K}(x)(G|-).$$

\noindent The corresponding formula then holds  for the coefficients of $x^{m}$
in the power series expansion of $F_{K}(x).$ \vspace{2mm}

\noindent {\bf Lemma 10.}  Let $F_{K}(x) = f_{K}(e^{x})$ denote the power series
resulting from substitution of $e^{x}$ for the variable $A$ in the Laurent
polynomial $f_{K}(A)$ (defined in section 2).  Write this power series in the
form

$$F_{K}(x) = \sum_{m=0}^{\infty}  v_{m}(K) x^{m}.$$

\noindent Then the numerical invariants $v_{k}(K)$ are of finite graphical type
$k$. \vspace{3mm}

\noindent {\bf Proof.}  Recall from section 5 that

$$f_{K_{+}} = -A^{-2} f_{K_{0}} - A^{-4} f_{K_{\infty}}$$ $$f_{K_{-}} = -A^{+2}
f_{K_{0}} - A^{+4} f_{K_{\infty}}$$

\noindent It follows that

$$F_{K_{*}}=  F_{K_{+}} - F_{K_{-}}$$

$$= f_{K_{+}}(e^{x}) - f_{K_{-}}(e^{x})$$

\noindent is divisible by $x.$  Thus if $G$ has $m$ nodes then $F_{G}$ is
divisible by $x^{m}.$  This implies, for any $G$, that $v_{k}(G) = 0$ if $k < m =
N(G).$  This is exactly the statement that $v_{k}$ is of finite graphical type
$k.$// \vspace{3mm}

\noindent {\bf Proposition 11.}  Let $G$ be a graph with n vertices so that
$N(G)=n$, configured as a virtual diagram in the plane.  Let $(G|+)$ denote  the
diagram  $G$ with a specific crossing of positive type and $(G|-)$ the diagram
identical to $G$ except that the crossing has been switched to one of negative
type. Let $(G|*)$ denote the result of replacing this crossing by a graphical
vertex.  Let $v$ be a Vassiliev invariant of type $n = N(G).$  Then $v(G|+) =
v(G|-).$  Thus a Vassiliev invariant of type $n$ is independent of the settings
of the crossings (plus or minus) in a diagram for $G.$ \vspace{3mm}

\noindent {\bf Proof.}  $v(G|+) - v(G|-) = v(G|*)$ by the definition of a
Vassiliev invariant.  But $v(G|*) = 0$ since $(G|*)$  has $(n+1)$ vertices and
$v$ is of type $n.$  This completes the proof.// \vspace{3mm}

\noindent {\bf Corollary 12. } If $G$ and $v$ are as in Proposition 11 , and $G$
is a classical diagram (free of virtual crossings) then $v(G)$ does not depend
upon the classical embedding of $G$ in $R^{3}$ that is indicated by the diagram.
\vspace{3mm}

\noindent {\bf Proof.} This follows directly from the switching independence
shown in Proposition 11. // \vspace{3mm}

For virtual Vassiliev invariants one should not expect the analog of this
corollary to hold, but in fact it does hold for the virtual Vassiliev invariants
induced from $f_{K}(A).$  That is, we shall show that the Vassiliev invariants
$v_{n}(G)$ in the series $F_{K}(x)$ depend only on the chord diagram associated
with $G$ when $G$ is a virtual diagram with $n$ graphical nodes.  This is the
subject of the following subsection. \vspace{3mm}

\subsection{The Vassiliev Invariants Induced by the Jones-Polynomial}

We shall use the Vassiliev invariants that arise from the bracket polynomial and
the $f$-polynomial. This is equivalent to using the Vassiliev invariants that
arise from the Jones polynomial.  Let f(A) be any Laurent polynomial with
coefficients as in the formula below

$$f(A) = c_{1}A^{d_{1}}  +  c_{2}A^{d_{2}}  + ... + c_{k}A^{d_{k}},$$

\noindent where  the degrees are integers arranged so that

$$d_{1} \leq d_{2} \leq ... d_{k}.$$

\noindent Then

$$f(e^{x}) = c_{1}e^{xd_{1}}  +  c_{2}e^{xd_{2}}  + ... + c_{k}e^{xd_{k}}$$

$$= \sum_{n=0}^{\infty} ( c_{1}d_{1}^{n}  +  c_{2}d_{2}^{n}  + ... +
c_{k}d_{k}^{n})x^{n}/n!.$$

\noindent Thus if

$$F(x) = f(e^{x}) = \sum_{n=0}^{\infty} v_{n} x^{n},$$

\noindent then

$$v_{n} = (c_{1}d_{1}^{n}  +  c_{2}d_{2}^{n}  + ... + c_{k}d_{k}^{n})/n!.$$

\noindent This gives a direct formula for the Vassiliev invariants $v_{n}$
associated with $f$. \vspace{3mm}

In particular, this gives us a direct method to read off the Vassiliev invariants
associated with a given evaluation of the normalized bracket polynomial $f_{K}$. 
The invariant $v_{n}(G)$ is determined by the coefficients of $f_{G}(A)$ and the
exponents of $A$ in this Laurent polynomial. \vspace{3mm}

\noindent {\bf Notational Discussion.}  Let $v_{n}(K)$ denote the $n$-th
Vassiliev invariant induced from $f_{K}(A)$ as described in this section. Let
$G_{*}$, $G_{+}$, $G_{-}$ denote a triple of (virtual) graph diagrams that differ
at the site of one rigid vertex (denoted $*$) by replacement by either a positive
crossing (denoted $+$) or a negative crossing (denoted $-$). Let $G_{0}$ and
$G_{\infty}$ denote the oriented and unoriented smoothings of this crossing. 
Note that since we can speak of the evaluation of $f_{G_{\infty}}(A)$, it follows
that $v_{n}$ is defined for diagrams with non-oriented smoothings - one just
evaluates the state sum in the usual way with single reverse-oriented loops
taking the usual loop value of $-A^{2} - A^{-2}.$ \vspace{3mm}

\noindent {\bf Theorem 13.}  With notation as above, the following recursion
formula holds for the Vassiliev invariants $v_{n}(G).$

$$v_{n}(G_{*}) = \sum_{k=0}^{n-1} c_{n,k}(v_{k}(G_{0}) +
2^{n-k}v_{k}(G_{\infty}))$$

\noindent where

$$(2^{n-k}(1+(-1)^{n-k+1})/(n-k)! = c_{n,k}.$$

\noindent The value of $v_{0}(K)$ on a virtual diagram without graphical nodes
depends only on the number of components in the diagram, and is independent of
the configuration of virtual crossings. Specifically,

$$v_{0}(K) = (-2)^{\mu(K) -1}$$

\noindent where  $\mu(K)$ denotes the number of link components in $K.$
\vspace{3mm}

\noindent {\bf Corollary 14.} The Vassiliev invariants $v_{n}(G)$ in the series
$F_{K}(x)$ depend only on the chord diagram associated with $G$ when $G$ is a
virtual diagram with $n$ graphical nodes. Hence the weight systems for the
invariants $v_{n}(G)$ do not depend upon virtual crossings. \vspace{3mm}

\noindent {\bf Proof of Corollary.}  This follows directly from Theorem 13 since
the recursion formula in that theorem computes $v_{n}(G)$ for a graph $G$ with n
nodes in terms of $v_{0}(K)$ for a collection of virtual knots $\{ K \}.$ Since
$v_{0}$ is independent of virtual crossings, so is  $v_{n}(G).$  This completes
the proof.// \vspace{3mm}

\noindent {\bf Proof of Theorem 13.}

Recall from section 5 that

$$f_{G_{+}} =-A^{-2} f_{G_{0}} - A^{-4} f_{G_{\infty}}$$

$$f_{G_{-}} = -A^{+2} f_{G_{0}} -A^{+4} f_{G_{\infty}}$$

\noindent Hence

$$f_{G_{*}} = (A^{2} - A^{-2})f_{G_{0}} + (A^{4} - A^{-4})f_{G_{\infty}}.$$

\noindent Now suppose that

$$f_{G_{0}} = \sum_{i} a_{i}A^{n_{i}}$$

\noindent and that

$$f_{G_{\infty}} = \sum_{i} b_{i}A^{m_{i}}.$$

\noindent Then

$$f_{G_{*}} = \sum_{i} a_{i}A^{n_{i} +2} -  a_{i}A^{n_{i} -2} + b_{i}A^{m_{i} +
4} - b_{i}A^{m_{i} - 4}.$$

\noindent Therefore

$$v_{n}(G_{*}) = (1/n!)\sum_{i} a_{i}((n_{i} +2)^{n} - (n_{i} -2)^{n}) +
b_{i}((m_{i} +4)^{n} - (m_{i} -4)^{n}).$$

\noindent The first part of Theorem 13 follows from this formula by a direct
application of the binomial theorem. For the second part, note that

$$v_{0}(K) = f_{K}(1) = (-1)^{w(K)}<K>(1).$$

\noindent For $A=1$ it is easy to see that the only effect of the matrix model of
section 6 on the bracket calculation is to multiply by $(-1)^{cv(K)}$ where
$cv(K)$ is the number of virtual crossings in $K.$ That is,

$$d<K>(1) = (-1)^{cv(K)}Z(K)(1) $$

\noindent where $Z(K)(1)$ is the matrix model of section 6 evaluated at $A=1.$
\vspace{3mm}

\noindent In this matrix model there is no difference (at $A=1$) between the
crossings and the virtual crossings. They are both algebraically crossed
Kronecker deltas. Consequently,  $Z(K) = Z(K')$ where $K'$ is the same diagram as
$K$ with all the virtual crossings replaced by (flat) classical crossings. Then 
it  follows from the regular isotopy invariance of  $Z$ that $Z(K')(1)
=(-2)^{\mu(K)} (-1)^{c(K')}$ where $c(K')$ is the total number of crossings in
$K'.$  Note that the value of $d = -2$ when $A=1.$ Hence

$$v_{0}(K) = (-1)^{w(K)} (-1)^{cv(K)} (-1)^{c(K')} (-2)^{\mu(K) -1}$$

\noindent Now we know that $w(K) + cv(K)$ is congruent modulo $2$ to $c(K').$
Therefore

$$v_{0}(K) = (-2)^{\mu(K) -1}.$$

\noindent This completes the proof of the Theorem.// \vspace{3mm}

%\begin{figure}[htbp] \vspace*{60mm} \special{pntg=F25.pntg} \vspace*{13pt}
%\begin{center} { Figure 25 --- $v_{2}$ dependence } \end{center} \end{figure}
%\vspace{3mm}

$$ \picill3inby2.5in(F25.EPSF) $$
\begin{center}
{\bf Figure 25 --- $v_{2}$ dependence }
\end{center}

\noindent {\bf Remark.} An example for this Theorem is that $v_{2}$ gives the
value $-48$ for both of the graphs shown in Figure 25.  Each graph has two
graphical nodes. One graph represents a virtual diagram inequivalent to any
embedding of the other. The invariants $v_{n}(G)$ themselves depend on virtual
crossings for the graphs with less than n nodes. In fact, in these intermediate
ranges there is a dependency on infinitely many virtual diagrams so that these
invariants are no longer as ``finite'' as the classical Vassiliev invariants.  In
\cite{GPV} there is formulated a more restrictive notion of finite type virtual
invariants. Our $v_{2}$ from the Jones polynomial is a first example of a finite
graphical type invariant that is outside the scheme proposed by Goussarov, Polyak
and Viro.  More work needs to be done to have a complete theory of virtual
Vassiliev invariants. \vspace{3mm}

\section{Discussion}

This completes our introduction to virtual knot theory. There is much that begs
for further investigation.  We leave the following topics for sequels to this
paper: the Alexander polynomial (there are a diversity of definitions that differ
on virtual knots), virtual braids, virtual 3-manifolds, Vassiliev invariants
induced from quantum link invariants, more general structure of Vassiliev
invariants. \vspace{2mm}

It should be remarked that the usual argument that induces Vassiliev invariants
from quantum link invariants produces virtual Vassiliev invariants from our
natural extension of quantum link invariants. Of course, we have to handle the
virtual framing for these cases as was discussed in section 6.  The matter of the
virtual framing needs further thought since introducing it means that we are no
longer just considering abstract Gauss codes. \vspace{2mm}

For general Vassiliev invariants it is worth comparing our results with those of
Goussarov, Polyak and Viro \cite{GPV}. The general notion of finite graphical
type given here and their notion of finite type suggest a unification not yet
fully perceived. \vspace{2mm}

Virtual braids is a subject very close to the ``welded braids" of Fenn, Rimanyi
and Rourke \cite{FRR}.  In fact, their welded braids are a quotient of the
category of virtual braids that are defined through our approach.  I am indebted
to Tom Imbo for pointing out this connection. This topic will be the subject of a
separate paper. \vspace{2mm}

This work began with an attempt to understand the Jones polynomial for classical
knots by generalising that category.  I hope that these considerations will lead
to deeper insight into the Jones polynomial and  its relationship with the
fundamental group and quandle of a classical knot. 

\end{document}